 \documentclass{amsart}    \usepackage[dvips]{hyperref}

\title{
	R\'{e}nyi dimension and Gaussian filtering
}
\author{
	Terry A. Loring
}
\address{
	Department of Mathematics and Statistics, University of New Mexico,
	Albuquerque, NM 87131, USA.
}
\email{
	loring@math.unm.edu
}
\thanks{
	This work was supported in part by DARPA Contract N00014-03-1-0900.
}

\theoremstyle{plain}
\newtheorem{thm}{Theorem}[section]
\newtheorem{lem}[thm]{Lemma}
 
\theoremstyle{definition}
\newtheorem{defn}[thm]{Definition}
\newtheorem{rem}[thm]{Remark}

\begin{document}

\begin{abstract}

Consider the partition function $S_{\mu}^{q}(\epsilon)$ associated
in theory of R\'{e}nyi dimension to a finite Borel measure $\mu$
on Euclidean d-space. This partion function $S_{\mu}^{q}(\epsilon)$
is the sum of the $q$-th powers of the measure applied to a partition
of $d$-space into $d$-cubes of width $\epsilon.$ We further Guerin's
investigation of the relation between this partition function and
the Lebesgue $L^{p}$ norm ($L^{q}$ norm) of the convolution of $\mu$
against an approximate identity of Gaussians. We prove a Lipschitz-type
esimate on the partition function. This bound on the partition function
leads to results regarding the computation of R\'{e}nyi dimension.  It
also shows that the partion function is of $O$-regular variation.

We find situtations where one can or cannot replace the partition
function by a discrete version. We discover that the slopes of the
least-square best fit linear approximations to the partion function
cannot always be used to calculate upper and lower R\'{e}nyi dimension.
\end{abstract}

\maketitle

\begin{center} (preprint version) \end{center}
\markright {R\'{e}nyi dimension and Gaussian filtering (preprint version)}
\markleft {R\'{e}nyi dimension and Gaussian filtering (preprint version)}

\section{Introduction}

The R\'{e}nyi dimensions of a finite Borel measure $\mu$ on $\mathbb{R}^{d}$
are derived from slopes of certain long secants of the log-log plot
of the function
\[
\epsilon\mapsto S_{\mu}^{q}(\epsilon)
=\sum_{\mathbf{k}\in\mathbb{Z}^{d}}\mu(\epsilon\mathbf{k}+\epsilon\mathbb{I})^{q},
\]
where
\[
\mathbb{I}=[0,1)\times[0,1)\times\cdots\times[0,1).
\]
There are exceptions for $q=0,1.$ In this paper we only address the
cases $0<q<1$ and $1<q<\infty.$ In this introduction we wish to
avoid convergence issues, so let us also assume that $\mu$ has bounded
support.

For any $x_{0},$ we set
\[
D_{q}^{\pm}(\mu)=
\lim_{x\rightarrow-\infty}\,_{\textrm{inf }}^{\textrm{sup}}
\frac{1}{q-1}\frac{\ln\left(S_{\mu}^{q}(e^{x})\right)
-\ln\left(S_{\mu}^{q}(e^{x_{0}})\right)}{x-x_{0}}.
\]
The constant terms are irrelevant, so this is usually written as
\[
D_{q}^{\pm}(\mu)=\lim_{x\rightarrow-\infty}\,_{\textrm{inf }}^{\textrm{sup}}
\frac{1}{q-1}\frac{\ln\left(S_{\mu}^{q}(e^{x})\right)}{x}
\]
or
\[
D_{q}^{\pm}(\mu)=\lim_{\epsilon\rightarrow0}\,_{\textrm{inf }}^{\textrm{sup}}
\frac{1}{q-1}\frac{\ln\left(\sum_{\mathbf{k}\in\mathbb{Z}^{d}}\mu(\epsilon\mathbf{k}
+\epsilon\mathbb{I})^{q}\right)}{\ln(\epsilon)}.
\]
Moving an exponent inside the log gives
\[
D_{q}^{\pm}(\mu)=\lim_{\epsilon\rightarrow0}\,_{\textrm{inf }}^{\textrm{sup}}
\frac{q}{q-1}\frac{\ln\left(\left\Vert \mathbf{k}\mapsto\mu(\epsilon\mathbf{k}
+\epsilon\mathbb{I})\right\Vert _{q}\right)}{\ln(\epsilon)}.
\]
This shows a relationship between convolution, $L^{p}$-norms and
R\'{e}nyi dimension, because 
\[
(\chi|_{(-\epsilon\mathbb{I})}\ast\mu)(\mathbf{x})=\mu(\mathbf{x}+\epsilon\mathbb{I}).
\]
Here we have used $\chi|_{(-\epsilon\mathbb{I})}$ to denote the characteristic
function of $-\epsilon\mathbb{I}.$

Gu\'{e}rin (\cite{Guerin}) showed a more general relation between
convolutions, $L^{p}$-norms and R\'{e}nyi dimensions. He showed
that for many choices of a scalar-valued function $g$ on $\mathbb{R}^{d},$
if
\[
1<q<\infty,
\]
and if we set
\[
g_{\epsilon}(\textrm{x})=\epsilon^{-d}g(\epsilon^{-1}\mathbf{x}),
\]
then
\[
D_{q}^{\pm}(\mu)=\lim_{\epsilon\rightarrow0}\,_{\textrm{inf }}^{\textrm{sup}}
\frac{1}{q-1}\frac{\ln\left(\epsilon^{d(q-1)}\left\Vert g_{\epsilon}\ast\mu
\right\Vert _{q}^{q}\right)}{\ln(\epsilon)},
\]
or
\begin{equation}
D_{q}^{\pm}(\mu)=d+\lim_{\epsilon\rightarrow0}\,_{\textrm{inf }}^{\textrm{sup}}
\frac{q}{q-1}\frac{\ln\left(\left\Vert 
g_{\epsilon}\ast\mu\right\Vert _{q}\right)}{\ln(\epsilon)}.
\label{eq:  Guerin's formula}
\end{equation}
Gu\'{e}rin allowed $g$ from a large class of complex-valued, rapidly
decreasing functions. 

A technical improvement on Gu\'{e}rin's result is given in Section
\ref{sec:Norms-after-Convolution}, with additional restrictions on
$g$ but allowing $0<q<1.$ For a given $\mu,$ and a ``nice'' function
$g\geq0,$ we establish a uniform bound on difference
\[
\ln\left(\epsilon^{d(q-1)}\left\Vert g_{\epsilon}\ast\mu\right\Vert _{q}^{q}\right)
-\ln\left(S_{\mu}^{q}(\epsilon)\right).
\]
This estimate allows us to analyze sequences
\[
\frac{\ln\left(S_{\mu}^{q}(\epsilon_{n})\right)}{\ln(\epsilon_{n})}
\]
by looking instead at
\[
\frac{\ln\left(\left\Vert 
g_{\epsilon_{n}}\ast\mu\right\Vert _{q}\right)}{\ln(\epsilon_{n})}.
\]
This will be advantageous if we choose $g_{\epsilon}$ properly.

Most importantly, we wish to let $g$ be a standard Gaussian on $\mathbb{R}^{d}.$
As we have the convention
\[
g_{\epsilon}(\textrm{x})=\epsilon^{-d}g(\epsilon^{-1}\mathbf{x}),
\]
the semigroup rule ends up as
\[
g_{\epsilon}\ast g_{\eta}=g_{\sqrt{\epsilon^{2}+\eta^{2}}}.
\]

We find that $\left\Vert g_{\epsilon}\ast\mu\right\Vert _{q}$
gives us information not apparent in the sums $S_{\mu}^{q}(\epsilon).$ 
Specifically, we find constants $A$ and $B$ so that
\[
\left|\ln\left(S_{\mu}^{q}(e^{x})\right)-\ln\left(S_{\mu}^{q}(e^{y})\right)\right|
\leq A+B|x-y|.
\]
It follows that $S_{\mu}^{q}(x^{-1})$ is of $O$-regular
variation.

Generalizing a result of Riedi (\cite{Riedi}), we show that
\[
\frac{1}{q-1}\frac{\ln\left(S_{\mu}^{q}(\epsilon_{n})\right)}{\ln(\epsilon_{n})}
\]
can be used to calculate $D_{q}^{\pm}(\mu),$ even if $\epsilon_{n}$
converges to zero somewhat faster than geometrically. The specific
requirement is that
\[
\lim_{n\rightarrow\infty}\frac{\ln(\epsilon_{n+1})}{\ln(\epsilon_{n})}>1.
\]
 Examples are exhibited that shows that this result is in some sense
the best possible.

We also give some new estimates on the R\'{e}nyi dimensions of a
convolution $\mu\ast\nu$ in terms of the R\'{e}nyi dimensions of
$\mu$ and $\nu.$

In Sections \ref{sec:Best-Fit-Slopes} and \ref{sec:Best-Fit-Slopes II}
we consider some of the changes that occur when one replaces
\[
\frac{\ln\left(S_{\mu}^{q}(e^{x})\right)
-\ln\left(S_{\mu}^{q}(e^{x_{0}})\right)}{x-x_{0}}
\]
by the slope of a least-squares best fit line over $[x,x_{0}]$ to
the function
\[
t\mapsto\ln\left(S_{\mu}^{q}(e^{t})\right).
\]
We exhibit an example where these least-squares slopes do not determine
the upper R\'{e}nyi dimension.

The examples we give have features that occur only on scales that
grow doubly exponentially. In the final section we suggest an alteration
of R\'{e}nyi that better detects the aberrant nature of these examples.

\section{Norms after convolution \label{sec:Norms-after-Convolution} }

Here follows our main technical result. Our initial interest here
was in the context of image analysis, were convolution by scaled Gaussians
is common. For example, see \cite{BurtAdelson}.

For $1<q<\infty,$ we have an easy finite bound on the partition function
\[
S_{\mu}^{q}(\epsilon)
=\sum_{\mathbf{k}\in\mathbb{Z}^{d}}\mu(\epsilon\mathbf{k}+\epsilon\mathbb{I})^{q},
\]
specifically
\begin{eqnarray*}
S_{\mu}^{q}(\epsilon) 
 & \leq & \left(\sum_{\mathbf{k}\in\mathbb{Z}^{d}}
 \mu(\epsilon\mathbf{k}+\epsilon\mathbb{I})\right)^{q}\\
 & = & \mu(\mathbb{R}^{d})^{q}\\
 & < & \infty.
 \end{eqnarray*}
For $0<q<1,$ it is possible to have $S_{\mu}^{q}(\epsilon)=\infty.$

\begin{defn}
A finite Borel measure $\mu$ on $\mathbb{R}^{d}$ is $q$\emph{-finite}
if $S_{\mu}^{q}(\epsilon)<\infty.$ Notice that if $\mu$ has bounded
support then $\mu$ is $q$-finite for all $0<q<1.$
\end{defn}

Barbaroux, Germinet, and Tcheremchantsev have the following result
implicitly in \cite{BarGerTch}. 

\begin{lem}
Let $\mu$ be a finite Borel measure on $\mathbb{R}^{d}$. For any
$0<q<1,$ the following are equivalent\/{\rm :}
\begin{enumerate}
\item $\mu$ is $q$-finite\/{\rm ;}
\item there exists $\epsilon>0$ for which $S_{\mu}^{q}(\epsilon)<\infty\/{\rm ;}$
\item for all $\epsilon>0$ it is true that $S_{\mu}^{q}(\epsilon)<\infty.$
\end{enumerate}
\end{lem}

\begin{proof}
On page 992-3 of \cite{BarGerTch} it is shown that if $\mu$ is $q$-finite
then $S_{\mu}^{q}(\epsilon)<\infty$ for small $\epsilon.$ Also,
it is shown that if $\mu$ is not $q$-finite then $S_{\mu}^{q}(\epsilon)=\infty$
for small $\epsilon.$ A rescaling argument show that if $S_{\mu}^{q}(\epsilon)$
is ever finite then it is finite for all small $\epsilon$, while
if it is ever infinite, it is infinite for all small $\epsilon.$
Therefore, the partition function is either finite for all $\epsilon$
or infinite for all $\epsilon.$
\end{proof}
The proof of the following borrows from the methods in \cite{BarGerTch}.

\begin{lem}
\label{lemma: norm of gaussians}
Suppose $g$ is a real-valued Borel
measurable function on $\mathbb{R}^{d}$ that is nonnegative, bounded,
bounded away from zero in a neighborhood of $\mathbf{0},$ and rapidly
decreasing. Let
\[
g_{\epsilon}(\mathbf{x})=\epsilon^{-d}g(\epsilon^{-1}\mathbf{x}).\]
Suppose $\mu$ is a finite Borel measure on $\mathbb{R}^{d}.$ If
$1<q<\infty,$ or if $0<q<1$ and $\mu$ is $q$-finite, then there
exists a constant $1<C<\infty$ so that for all positive $\epsilon,
$\begin{equation}
C^{-1}\leq\frac{\epsilon^{d(q-1)}\left\Vert 
g_{\epsilon}\ast\mu\right\Vert _{q}^{q}}{S_{\mu}^{q}(\epsilon)}\leq C.
\label{eq: ratio}
\end{equation}
Here the $L^{q}$ norm is with respect to Lebesgue measure. Therefore
\[
D_{q}^{\pm}(\mu)=d+\lim_{\epsilon\rightarrow0}\,_{\textrm{inf }}^{\textrm{sup}}
\frac{q}{q-1}\frac{\ln\left(\left\Vert 
g_{\epsilon}\ast\mu\right\Vert _{q}\right)}{\ln(\epsilon)}.
\]
\end{lem}

\begin{proof}
We will use $m$ to denote Lebesgue measure, to keep it straight from
$\mu.$ Let us denote the open unit rectangle at the origin by $\mathbb{D},$
so
\[
\mathbb{D}=(-1,1)\times(-1,1)\times\cdots\times(-1,1).
\]
Recall $\mathbb{I}$ is the product of $d$ copies of $[0,1).$

Let us denote by $\mu^{(\epsilon)}$ the sequences over $\mathbb{Z}^{d}$
given by
\[
\mu_{\mathbf{n}}^{(\epsilon)}=\mu(\epsilon\mathbf{n}+\epsilon\mathbb{I}).
\]
Thus
\[
S_{\mu}^{q}(\epsilon)
=\sum_{\mathbf{k}\in\mathbb{Z}^{d}}\mu(\epsilon\mathbf{k}+\epsilon\mathbb{I})^{q}
=\left\Vert \mu^{(\epsilon)}\right\Vert _{q}^{q}
\]
(the norm here is on $l^{q}(\mathbb{Z}^{d})$).

An obvious rescaling reduces this theorem to the special case where
\[
\inf\{ g(\mathbf{x})\mid\mathbf{x}\in\mathbb{D}\}>0,
\]
so let us make this assumption. We compute
\begin{equation}
 \left\Vert g_{\epsilon}\ast\mu\right\Vert _{q}^{q}
 =  \epsilon^{-qd}\sum_{\mathbf{j}\in\mathbb{Z}^{d}}
 \int_{\epsilon\mathbf{j}+\epsilon\mathbb{I}}\left(\sum_{\mathbf{k}\in\mathbb{Z}^{d}}
 \int_{\epsilon\mathbf{k}+\epsilon\mathbb{I}}g(\epsilon^{-1}(\mathbf{x}-\mathbf{y}))
 \, d\mu(\mathbf{y})\right)^{q}\, dm(\mathbf{x}).
 \label{eq: cut space}
\end{equation}
If 
\[
\mathbf{x}\in\epsilon\mathbf{j}+\epsilon\mathbb{I}
\]
and
\[
\mathbf{y}\in\epsilon\mathbf{k}+\epsilon\mathbb{I}
\]
then
\[
\epsilon^{-1}(\mathbf{x}-\mathbf{y})\in(\mathbf{j}-\mathbf{k})+\mathbb{D}.
\]
Let us define $\gamma$ and $\Gamma$ as sequences over $\mathbb{Z}^{d},$
by
\[
\gamma_{\mathbf{n}}=\inf\{ g(\mathbf{x})\mid\mathbf{x}\in\mathbf{n}+\mathbb{D}\}
\]
and
\[
\Gamma_{\mathbf{n}}=\sup\{ g(\mathbf{x})\mid\mathbf{x}\in\mathbf{n}+\mathbb{D}\}.
\]
These give us bounds on the $g(\epsilon^{-1}(\mathbf{x}-\mathbf{y}))$
term inside integrals in (\ref{eq: cut space}).

For an upper bound, we get
\begin{eqnarray*}
\left\Vert g_{\epsilon}\ast\mu\right\Vert _{q}^{q} 
 & \leq & \epsilon^{-qd}\sum_{\mathbf{j}\in\mathbb{Z}^{d}}
 \left(\sum_{\mathbf{k}\in\mathbb{Z}^{d}}
 \Gamma_{\mathbf{j}-\mathbf{k}}\mu(\epsilon\mathbf{k}+\epsilon\mathbb{I})
 \right)^{q}\epsilon^{d}\\
 & = & \epsilon^{d(1-q)}\left\Vert 
 \Gamma\ast\mu^{(\epsilon)}\right\Vert _{q}^{q}.
\end{eqnarray*}
For a lower bound,
\[
\left\Vert g_{\epsilon}\ast\mu\right\Vert _{q}^{q}
\geq\epsilon^{d(1-q)}\left\Vert \gamma\ast\mu^{(\epsilon)}\right\Vert _{q}^{q}.
\]

First assume $1<q<\infty.$ The assumption that $g$ is nonzero on
$\mathbb{D}$ is here used to obtain
$\left\Vert \gamma\right\Vert _{q}>0.$
From the rapidly decreasing assumption we obtain
$\left\Vert \Gamma\right\Vert _{1}<\infty.$
Since,
\begin{eqnarray*}
\left\Vert g_{\epsilon}\ast\mu\right\Vert _{q} 
& \leq & \epsilon^{d\frac{(1-q)}{q}}\left\Vert 
\Gamma\ast\mu^{(\epsilon)}\right\Vert _{q}\\
 & \leq & \epsilon^{d\frac{(1-q)}{q}}\left\Vert 
 \Gamma\right\Vert _{1}\left\Vert \mu^{(\epsilon)}\right\Vert _{q}
 \end{eqnarray*}
and
\begin{eqnarray*}
\left\Vert g_{\epsilon}\ast\mu\right\Vert _{q} 
& \geq & \epsilon^{d\frac{(1-q)}{q}}\left\Vert 
\gamma\ast\mu^{(\epsilon)}\right\Vert _{q}\\
 & \geq & \epsilon^{d\frac{(1-q)}{q}}\left\Vert 
 \gamma\right\Vert _{q}\left\Vert \mu^{(\epsilon)}\right\Vert _{q},
 \end{eqnarray*}
we may take
\[
C=
\max\left(\left\Vert 
\Gamma\right\Vert _{1}^{q},\left\Vert \gamma\right\Vert _{q}^{-q}\right).
\]

Now assume $0<q<1.$ The assumptions on $g$ tell us 
$\left\Vert \gamma\right\Vert _{1}>0$
and 
$\left\Vert \Gamma\right\Vert _{q}<\infty.$ 
In this case we may take
\[
C=\max\left(
\left\Vert \Gamma\right\Vert _{q}^{q},\left\Vert \gamma\right\Vert _{1}^{-q}
\right).
\]

\end{proof}

\section{Bounding the partition function \label{sec:Bounding-the-Partition}}

In this section we use Lemma \ref{lemma: norm of gaussians} to establish
bounds on
\[
\ln\left(S_{\mu}^{q}(e^{x+h})\right)-\ln\left(S_{\mu}^{q}(e^{x})\right)
\]
that are of first order in $h$ and hold for all $x.$ Recall that the partition
function for $\mu$ is
\[
S_{\mu}^{q}(\epsilon)
=\sum_{\mathbf{k}\in\mathbb{Z}^{d}}\mu(\epsilon\mathbf{k}+\epsilon\mathbb{I})^{q}.
\]

Much of this section is familiar, but one conclusion seems novel:
$S_{\mu}^{q}(x^{-1})$ is almost decreasing. By this we mean
\[
\sup_{x<y}\frac{S_{\mu}^{q}(y^{-1})}{S_{\mu}^{q}(x^{-1})}<\infty.
\]
(See \cite{BingGoldTeu}.)

In this section, $g$ denotes a standard Gaussian, and for $\epsilon>0$
we set
\[
g_{\epsilon}(\mathbf{x})=\epsilon^{-d}g(\epsilon^{-1}\mathbf{x}).
\]
Let us also adopt the notation
\[
I_{\mu}^{q}(\epsilon)=\left\Vert g_{\epsilon}\ast\mu\right\Vert _{q}^{q}.
\]

\begin{lem}
\label{lem:I is decreasing}Suppose $\mu$ is a finite Borel measure.
If $1<q<\infty$ then
\[
\epsilon_{1}\leq\epsilon_{2}\implies\left\Vert 
g_{\epsilon_{2}}\ast\mu\right\Vert _{q}
\leq\left\Vert g_{\epsilon_{1}}\ast\mu\right\Vert _{q}.
\]
If $0<q<1,$ and $\mu$ is $q$-finite, then
\[
\epsilon_{1}\leq\epsilon_{2}
\implies
\left\Vert g_{\epsilon_{2}}\ast\mu\right\Vert _{q}
\geq\left\Vert g_{\epsilon_{1}}\ast\mu\right\Vert _{q}.
\]

\end{lem}

\begin{proof}
If $1<q<\infty$ then
\begin{eqnarray*}
\left\Vert g_{\eta}\ast\mu\right\Vert _{q}  
 & = & 
 \left\Vert g_{\sqrt{\eta^{2}-\epsilon^{2}}}\ast g_{\epsilon}\ast\mu\right\Vert _{q}\\
 & \leq &
 \left\Vert g_{\sqrt{\eta^{2}-\epsilon^{2}}}\right\Vert _{1}
 \left\Vert g_{\epsilon}\ast\mu\right\Vert _{q}\\
 & = & \left\Vert g_{\epsilon}\ast\mu\right\Vert _{q}.
 \end{eqnarray*}
For $0<q<1$ the middle step becomes
\[
\left\Vert g_{\sqrt{\eta^{2}-\epsilon^{2}}}\ast g_{\epsilon}\ast\mu\right\Vert _{q}
\geq\left\Vert g_{\sqrt{\eta^{2}-\epsilon^{2}}}\right\Vert _{1}\left\Vert 
g_{\epsilon}\ast\mu
\right\Vert _{q}.
\]

\end{proof}

\begin{lem}
\label{lem: exact jumps}
Suppose $\mu$ is a finite Borel measure
and that $n$ is a natural number. For $1<q<\infty,$ if $\epsilon>0$
then
\[
0 \leq
\ln\left(S_{\mu}^{q}(\epsilon)\right)-\ln\left(S_{\mu}^{q}(2^{-n}\epsilon)\right)
\leq n\ln\left(2^{d(q-1)}\right).
\]
For $0<q<1,$ if $\mu$ is $q$-finite and $\epsilon>0$ then
\[
n\ln\left(2^{d(q-1)}\right)
\leq\ln\left(S_{\mu}^{q}(\epsilon)\right)-\ln\left(S_{\mu}^{q}(2^{-n}\epsilon)\right)
\leq0.
\]

\end{lem}

\begin{proof}
Suppose $1<q<\infty.$ Given a disjoint union of Borel sets, 
\[
E_{1}\cup\cdots\cup E_{2^{d}}=F
\]
we have the estimates
\[
2^{d(1-q)}\mu(F)^{q}\leq\sum_{j=1}^{2^{d}}\mu(E_{j})^{q}\leq\mu(F)^{q}.
\]
Therefore
\[
\ln\left(S_{\mu}^{q}(\epsilon)\right)+\ln(2^{d(1-q)})
\leq\ln\left(S_{\mu}^{q}(2^{-1}\epsilon)\right)
\leq\ln\left(S_{\mu}^{q}(\epsilon)\right).
\]

For $0<q<1$ the inequalities are all easily reversed.
\end{proof}

\begin{lem}
\label{lem:small jumps}
Suppose $\mu$ is a finite Borel measure.
If $1<q<\infty$ then there is a constant $E$ so that
\[
|\ln(\epsilon)-\ln(\eta)|\leq\ln(2)
\]
implies 
\[
\left|\ln\left(S_{\mu}^{q}(\epsilon)\right)-\ln\left(S_{\mu}^{q}(\eta)\right)\right|
\leq E.
\]
 If $0<q<1,$ and $\mu$ is $q$-finite, then there is constant $E$
so that
\[
|\ln(\epsilon)-\ln(\eta)|\leq\ln(2)
\]
implies 
\[
\left|\ln\left(S_{\mu}^{q}(\epsilon)\right)-\ln\left(S_{\mu}^{q}(\eta)\right)\right|
\leq E.
\]

\end{lem}

\begin{proof}
We may assume
\[
2^{-1}\epsilon\leq\eta\leq\epsilon.
\]
By Theorem \ref{lemma: norm of gaussians}, there is a constant $C$
so that for all $\rho>0,$
\[
\left|d(q-1)\ln(\rho)
+\ln\left(I_{\mu}^{q}(\rho)\right)
-\ln\left(S_{\mu}^{q}(\rho)\right)\right|
\leq D.
\]
If we set
\[
Q=d(q-1)\ln\left(2^{-\frac{1}{2}}\epsilon\right)
\]
and
\[
D=d|q-1|\ln\left(2^{\frac{1}{2}}\right)
\]
then 
\[
2^{-1}\epsilon\leq\rho\leq\epsilon
\implies
\left|Q+\ln\left(I_{\mu}^{q}(\rho)\right)-\ln\left(S_{\mu}^{q}(\rho)\right)\right|
\leq D.
\]

Now assume $1<q<\infty.$ Lemma \ref{lem:I is decreasing} gives us
\[
\ln\left(I_{\mu}^{q}(\epsilon)\right)\leq\ln\left(I_{\mu}^{q}(\eta)\right)
\leq\ln\left(I_{\mu}^{q}(2^{-1}\epsilon)\right)
\]
and Lemma \ref{lem: exact jumps} gives us
\[
\ln\left(S_{\mu}^{q}(2^{-1}\epsilon)\right)
\leq\ln\left(S_{\mu}^{q}(\epsilon)\right).
\]
We put this information together as
\begin{eqnarray*}
\ln\left(S_{\mu}^{q}(\epsilon)\right)-2D_{1} 
& \leq & Q+\ln\left(I_{\mu}^{q}(\epsilon)\right)-D\\
 & \leq & Q+\ln\left(I_{\mu}^{q}(\eta)\right)-D\\
 & \leq & \ln\left(S_{\mu}^{q}(\eta)\right)\\
 & \leq & Q+\ln\left(I_{\mu}^{q}(\eta)\right)+D\\
 & \leq & Q+\ln\left(I_{\mu}^{q}(2^{-1}\epsilon)\right)+D\\
 & \leq & \ln\left(S_{\mu}^{q}(2^{-1}\epsilon)\right)+2D\\
 & \leq & \ln\left(S_{\mu}^{q}(\epsilon)\right)+2D
 \end{eqnarray*}

For $0<q<1$ the proof is similar.
\end{proof}

\begin{thm}
\label{thm: order 1 diff}
Suppose $\mu$ is a finite Borel measure
on $\mathbb{R}^{d}.$ Let $B=d(q-1).$ For $1<q<\infty,$ there is
constant $A$ so that
\[
0<\eta<\epsilon
\]
implies 
\[
-A\leq\ln\left(S_{\mu}^{q}(\epsilon)\right)-\ln\left(S_{\mu}^{q}(\eta)\right)
\leq A+B(\ln(\epsilon)-\ln(\eta)).
\]
 For $0<q<1,$ if $\mu$ is $q$-finite then there is constant $A$
so that
\[
0<\eta<\epsilon
\]
implies 
\[
-A+B(\ln(\epsilon)-\ln(\eta))
\leq\ln\left(S_{\mu}^{q}(\epsilon)\right)-\ln\left(S_{\mu}^{q}(\eta)\right)
\leq A.
\]

\end{thm}

\begin{proof}
Assume first that $1<q<\infty.$ For some natural number $n,$
\[
2^{-n-1}\epsilon\leq\eta\leq2^{-n}\epsilon.
\]
This means
\[
n\leq\frac{1}{\ln(2)}(\ln(\epsilon)-\ln(\eta))+1
\]
and so by the last three lemmas,
\begin{eqnarray*}
 &  & \ln\left(S_{\mu}^{q}(\epsilon)\right)-B(\ln(\epsilon)-\ln(\eta))-B\ln(2)-E\\
 & \leq & \ln\left(S_{\mu}^{q}(\epsilon)\right)-nB\ln(2)-E\\
 & \leq & \ln\left(S_{\mu}^{q}(2^{-n}\epsilon)\right)-E\\
 & \leq & S_{\mu}^{q}(\eta)\\
 & \leq & \ln\left(S_{\mu}^{q}(2^{-n}\epsilon)\right)+E\\
 & \leq & \ln\left(S_{\mu}^{q}(\epsilon)\right)+E\\
 & \leq & \ln\left(S_{\mu}^{q}(\epsilon)\right)+B\ln(2)+E
 \end{eqnarray*}
so we can set
\[
A=B\ln(2)+E.
\]

For $0<q<1$ the proof is similar.
\end{proof}

\section{Application to discrete limits}

Riedi \cite{RiediDisseration,Riedi} shows that 
\begin{equation}
D_{p}^{\pm}(\mu)
=\lim_{n\rightarrow\infty}\,_{\textrm{inf }}^{\textrm{sup}}
\frac{1}{q-1}\frac{\ln\left(S_{\mu}^{q}(\epsilon_{n})\right)}{\ln(\epsilon_{n})}
\label{eq: seq ok}
\end{equation}
for $\epsilon_{n}\searrow0,$ so long as
\[
\limsup_{n\rightarrow\infty}(\ln(\epsilon_{n})-\ln(\epsilon_{n+1}))<\infty.
\]
 Indeed, he works with all $p\in\mathbb{R},$ and shows that grids
other than 
\[
\{\epsilon\mathbf{k}+\epsilon\mathbb{I}\mid\mathbf{k}\in\mathbb{Z}^{d}\}
\]
can be used. What concerns us here is that Riedi showed that (\ref{eq: seq ok})
is valid for a geometric series. We can go further, to allow sequences
such a $\epsilon_{n}=b^{-n^{r}}.$

\begin{lem}
\label{lemma:  bounded adjustments give equivalent seqeunces}
Suppose
$\mu$ is a finite Borel measure on $\mathbb{R}^{d},$ and $q\neq1$
is a positive number. If $0<q<1$ then also suppose 
$\sum\mu(\mathbf{k}+\mathbb{I})^{q}$
is finite. If $\epsilon_{n}\searrow0$ and $\eta_{n}\searrow0$ with
\[
\lim_{n\rightarrow\infty}\frac{\ln(\eta_{n})}{\ln(\epsilon_{n})}=1,
\]
then
\[
\lim_{n\rightarrow\infty}
\left|\frac{\ln\left(S_{\mu}^{q}(\eta_{n})\right)}{\ln(\eta_{n})}
-\frac{\ln\left(S_{\mu}^{q}(\epsilon_{n})\right)}{\ln(\epsilon_{n})}\right|
=0
\]

\end{lem}

\begin{proof}
Let $A$ and $B$ be the constants from Theorem \ref{thm: order 1 diff},
and let
\[
A_{0}=A+\left|\ln(S_{\mu}^{p}(1)\right|
\]
so that
\[
\left|\ln\left(S_{\mu}^{q}(\epsilon)\right)\right|\leq A_{0}+C|\ln(\epsilon)|
\]
for all $\epsilon.$ We may assume $\epsilon_{n}\leq e^{-1}$ and
$\eta_{n}\leq e^{-1},$ in which case
\begin{eqnarray*}
 &  & 
 \left|\frac{\ln\left(S_{\mu}^{q}(\eta_{n})\right)}{\ln(\eta_{n})}
 -\frac{\ln\left(S_{\mu}^{q}(\epsilon_{n})\right)}{\ln(\epsilon_{n})}\right|\\
 & \leq & \frac{\left|\ln\left(S_{\mu}^{q}(\eta_{n})\right)
 -\ln\left(S_{\mu}^{q}(\epsilon_{n})\right)\right|}{-\ln(\eta_{n})}
 +\frac{\left|\ln\left(S_{\mu}^{q}(\epsilon_{n})\right)\right|}
 {-\ln(\epsilon_{n})}\left|\frac{\ln(\epsilon_{n})}{\ln(\eta_{n})}-1\right|\\
 & \leq & \frac{A+C|\ln(\eta_{n})-\ln(\epsilon_{n})|}{-\ln(\eta_{n})}
 +\frac{A_{0}-C\ln(\epsilon_{n})}{-\ln(\epsilon_{n})}
 \left|\frac{\ln(\epsilon_{n})}{\ln(\eta_{n})}-1\right|\\
 & \leq & \frac{-A}{\ln(\eta_{n})}+
 C\left|1 -\frac{\ln(\epsilon_{n})}{\ln(\eta_{n})}\right|
 +\left(A_{0}+C\right)\left|\frac{\ln(\epsilon_{n})}{\ln(\eta_{n})}-1\right|\\
 & \rightarrow & 0.
 \end{eqnarray*}

\end{proof}

\begin{thm}
\label{thm: slow sequences get accumulation points}Suppose $\mu$
is a finite Borel measure on $\mathbb{R}^{d}$ and $q\neq1$ is a
positive number. If $0<q<1$ then also suppose $\mu$ is $q$-finite.
If $\epsilon_{n}\searrow0$ and
\[
\lim_{n\rightarrow\infty}\frac{\ln(\epsilon_{n+1})}{\ln(\epsilon_{n})}=1,
\]
then
\[
D_{q}^{\pm}(\mu)
=\frac{1}{q-1}\lim_{n\rightarrow\infty}\,_{\textrm{inf }}^{\textrm{sup}}
\frac{\ln\left(S_{\mu}^{q}(\epsilon_{n})\right)}{\ln(\epsilon_{n})}.
\]
Moreover, the sequence
\begin{equation}
n\mapsto\frac{\ln\left(S_{\mu}^{q}(\epsilon_{n})\right)}{\ln(\epsilon_{n})}
\label{eq: n sequence}
\end{equation}
and the net
\begin{equation}
\epsilon\mapsto\frac{\ln\left(S_{\mu}^{q}(\epsilon)\right)}{\ln(\epsilon)}
\label{eq: epsilon net}
\end{equation}
have the same accumulation points in $[-\infty,\infty].$
\end{thm}

\begin{proof}
Suppose $x$ is an accumulation point of the net in (\ref{eq: epsilon net}).
This means there is a sequence $\eta_{n}\searrow0$ so that
\[
\lim_{n\rightarrow\infty}
\frac{\ln\left(S_{\mu}^{q}(\eta_{n})\right)}{\ln(\eta_{n})}=x.
\]
Let the sequence $k_{n}$ be defined so that
\[
\epsilon_{k_{n}+1}\leq\eta_{n}\leq\epsilon_{k_{n}}.
\]
Of course, $k_{n}\searrow0$ and 
\[
\lim_{n\rightarrow\infty}
\frac{\ln\left(\eta_{n}\right)}{\ln\left(\epsilon_{k_{n}}\right)}
=1.
\]
The last lemma tells us
\[
\lim_{n\rightarrow\infty}
\frac{\ln\left(S_{\mu}^{q}(\epsilon_{k_{n}})\right)}{\ln(\epsilon_{k_{n}})}
=\lim_{n\rightarrow\infty}
\frac{\ln\left(S_{\mu}^{q}(\eta_{n})\right)}{\ln(\eta_{n})}
=x.
\]
Thus $x$ is also an accumulation point of the sequence in (\ref{eq: n sequence}).
\end{proof}

\section{Examples\label{sec:Examples}}

It is possible to use simple recursive definitions to create a Borel
measure $\mu,$ with support in $[0,1],$ so that the partition function
\[
S_{\mu}^{q}(\epsilon)=
\sum_{k\in\mathbb{Z}}\mu(\epsilon k+\epsilon\mathbb{I})^{q}
\]
behaves almost any way we would like. However, we do need to respect
Lemma \ref{lem: exact jumps}.

\begin{lem}
\label{lem:a_n example}Suppose 
\[
0\leq a_{n}\leq1
\]
for $n\geq1.$  If $q\neq1$ is a positive real number, there is a
Borel probability measure $\mu$ on $[0,1]$ for which
\[
\ln\left(S_{\mu}^{q}(2^{-n})\right)=(1-q)\ln(2)\sum_{j=1}^{n}a_{n}.
\]
Moreover, the net
\[
\epsilon\mapsto\frac{1}{q-1}
\frac{\ln\left(S_{\mu}^{q}(\epsilon)\right)}{\ln(\epsilon)}
\]
has the same accumulation points as the sequence
\[
n\mapsto\frac{1}{n}\sum_{j=1}^{n}a_{j}.
\]

\end{lem}

\begin{proof}
First let's define $F,$ the cumulative distribution function. We
start with 
\[
F(0)=0\quad\textrm{and}\quad F(1)=1.
\]
Choose any $\omega_{n}$ in $\left[0,\frac{1}{2}\right]$ so that
\[
\ln\left(\omega_{n}^{q}+(1-\omega_{n})^{q}\right)=(1-q)\ln(2)a_{n}.
\]

Define $F$ on the dyadic rationals between $0$ and $1$ by
\[
F\left(\frac{2k+1}{2^{n+1}}\right)
=\omega_{n}F\left(\frac{k}{2^{n}}\right)+
(1-\omega_{n})F\left(\frac{k+1}{2^{n}}\right).
\]
Thus $F$ is nondecreasing on the dyadic rationals; set it to $0$
on dyadics less than $0$ and to one on dyadics greater than $1.$
For any $n$ and any $k$ with $0<k<2^{-n}$ we have
\begin{eqnarray*}
 &  &
 F\left(\frac{k}{2^{n}}\right)-F\left(\frac{k}{2^{n}}-\frac{1}{2^{n+1}}\right)\\
 & = & F\left(\frac{k}{2^{n}}\right)
 -\left(\omega_{n}F\left(\frac{k-1}{2^{n}}\right)
 +(1-\omega_{n})F\left(\frac{k}{2^{n}}\right)\right)\\
 & = & \omega_{n}\left(F\left(\frac{k}{2^{n}}\right)
 -F\left(\frac{k}{2^{n}}-\frac{1}{2^{n}}\right)\right).
 \end{eqnarray*}
Since we choose $\omega_{n}$ to the left of $\frac{1}{2}$, we have
\[
F\left(\frac{k}{2^{n}}\right)
-F\left(\frac{k}{2^{n}}
-\frac{1}{2^{n+1}}\right)\leq\frac{1}{2}\left(F\left(\frac{k}{2^{n}}\right)
-F\left(\frac{k}{2^{n}}-\frac{1}{2^{n}}\right)\right).
\]
Therefore, if $m\geq1,$
\[
F\left(\frac{k}{2^{n}}\right)-F\left(\frac{k}{2^{n}}
-\frac{1}{2^{n+m}}\right)
\leq\left(\frac{1}{2}\right)^{m}\left(F\left(\frac{k}{2^{n}}\right)
-F\left(\frac{k}{2^{n}}-\frac{1}{2^{n}}\right)\right).
\]
Since $F$ is bounded between $0$ and $1,$ we have for any dyadic
rational $r,$ 
\[
\frac{k}{2^{n}}-\frac{1}{2^{n+m}}\leq r<\frac{k}{2^{n}}
\implies
F\left(\frac{k}{2^{n}}\right)-F\left(r\right)\leq\left(\frac{1}{2}\right)^{m}.
\]
Since $F$ is nondecreasing, this says
\[
F(s)=\sup_{s>r\in\mathbb{Z}\left[\frac{1}{2}\right]}F(r)
\]
for all $s$ in $\mathbb{Z}[\frac{1}{2}].$

Let us extend $F$ to $\mathbb{R}$ by 
\[
F(x)=\sup_{x>r\in\mathbb{Z}\left[\frac{1}{2}\right]}F(r).
\]
It is routine to verify that $F$ is left continuous and nondecreasing.

Since $F$ is nondecreasing and left continuous, we have an associated
measure $\mu$ which satisfies
\[
\mu([a,b))=F(b)-F(a).
\]
By design,
\begin{eqnarray*}
 &  & \mu\left(\left[\frac{2k}{2^{n+1}},\frac{2k+1}{2^{n+1}}\right)\right)\\
 & = & F\left(\frac{2k+1}{2^{n+1}}\right)
 -F\left(\frac{2k}{2^{n+1}}\right)\\
 & = & \omega_{n}F\left(\frac{k}{2^{n}}\right)
 +(1-\omega_{n})F\left(\frac{k+1}{2^{n}}\right)-F\left(\frac{k}{2^{n}}\right)\\
 & = &
 (1-\omega_{n})\mu\left(\left[\frac{k}{2^{n}},\frac{k+1}{2^{n}}\right)\right)
 \end{eqnarray*}
and
\begin{eqnarray*}
 &  & \mu\left(\left[\frac{2k+1}{2^{n+1}},\frac{2k+2}{2^{n+1}}\right)\right)\\
 & = & F\left(\frac{2k+2}{2^{n+1}}\right)-F\left(\frac{2k+1}{2^{n+1}}\right)\\
 & = & F\left(\frac{k+1}{2^{n}}\right)
 -\left(\omega_{n}F\left(\frac{k}{2^{n}}\right)
 +(1-\omega_{n})F\left(\frac{k+1}{2^{n}}\right)\right)\\
 & = & \omega_{n}\mu\left(\left[\frac{k}{2^{n}},\frac{k+1}{2^{n}}\right)\right).
 \end{eqnarray*}

This tells us
\[
S_{\mu}^{q}(2^{-n-1})=(\omega_{n}^{q}+(1-\omega_{n})^{q})S_{\mu}^{q}(2^{-n})
\]
and so also
\[
\ln\left(S_{\mu}^{q}(2^{-n-1})\right)
=(1-q)\ln(2)a_{n}+\ln\left(S_{\mu}^{q}(2^{-n})\right).
\]
Since 
\[
S_{\mu}^{q}(1)=\|\mu\|_{1}=1,
\]
induction gives us
\[
\ln\left(S_{\mu}^{q}(2^{-n})\right)=(1-q)\ln(2)\sum_{j=1}^{n}a_{j}.
\]

With $\mu$ as constructed from $a_{n}$ and $q$ as indicated, Theorem
\ref{thm: slow sequences get accumulation points} Theorem applies
to give the final statement in the lemma.
\end{proof}

\begin{thm}
Suppose $0<q<1$ or $1<q<\infty,$ and $\epsilon_{n}\searrow0$ with
\[
\liminf_{n\rightarrow\infty}\frac{\ln(\epsilon_{n+1})}{\ln(\epsilon_{n})}>1.
\]
Then there is a Borel probability measure on $[0,1]$ so that
\[
\lim_{n\rightarrow\infty}
\frac{\ln\left(S_{\mu}^{q}(\epsilon_{n})\right)}{\ln(\epsilon_{n})}
\]
exists but
\[
\lim_{\epsilon\rightarrow0}
\frac{\ln\left(S_{\mu}^{q}(\epsilon)\right)}{\ln(\epsilon)}
\]
does not.
\end{thm}

\begin{proof}
By Lemma \ref{lemma:  bounded adjustments give equivalent seqeunces},
it suffices to prove this in the special case where
\[
\epsilon_{n}=4^{-k_{n}}
\]
for some $k_{n}\in\mathbb{N}.$ The hypothesis on the $\epsilon_{n}$
translates to the assumption that $k_{n}$ is nondecreasing, with
limit $\infty,$ and that
\[
\liminf_{n\rightarrow\infty}\frac{k_{n+1}}{k_{n}}>1.
\]

Select a subsequence $k_{n_{l}}$ and some $R>1$ so that
\[
k_{n_{l}+1}\geq Rk_{n_{l}}
\]
for all $l.$ Define $a_{j}=\frac{1}{2}$ for all $j$ except
\[
2k_{n_{l}}<j\leq(k_{n_{l}}+k_{n_{l}+1})\implies a_{j}=0
\]
and
\[
(k_{n_{l}}+k_{n_{l}+1})<j\leq2k_{n_{l}+1}\implies a_{j}=1
\]
With $\mu$ as defined above,
\[
\frac{1}{q-1}
\frac{\ln\left(S_{\mu}^{p}(\epsilon_{n})\right)}{\ln(\epsilon_{n})}
=\frac{1}{2k_{n}}\sum_{j=1}^{2k_{n}}a_{j}=\frac{1}{2}.
\]

If 
\[
\eta_{n}=2^{-(k_{n_{l}}+k_{n_{l}+1})}
\]
then
\begin{eqnarray*}
\frac{1}{q-1}\frac{\ln\left(S_{\mu}^{q}(\eta_{n})\right)}{\ln(\eta_{n})} 
& = & \frac{2k_{n_{l}}}{k_{n_{l}}+k_{n_{l}+1}}\frac{1}{2}
+\frac{k_{n_{l}+1}-k_{n_{l}}}{k_{n_{l}}+k_{n_{l}+1}}0\\
 & = & \frac{1}{1+\frac{k_{n_{l}+1}}{k_{n_{l}}}}\\
 &  & \geq\frac{1}{1+R}.
 \end{eqnarray*}
Thus 
\begin{eqnarray*}
\liminf_{\epsilon\rightarrow0}
\frac{\ln\left(S_{\mu}^{q}(\epsilon)\right)}{\ln(\epsilon)} 
& \leq &
 \liminf_{n\rightarrow\infty}
 \frac{\ln\left(S_{\mu}^{q}(\epsilon_{n})\right)}{\ln(\epsilon_{n})}\\
 & = &
 \frac{1}{2}
\end{eqnarray*}
an\begin{eqnarray*}
\limsup_{\epsilon\rightarrow0}\frac{\ln\left(S_{\mu}^{q}(\epsilon)\right)}{\ln(\epsilon)}
& \geq & 
\limsup_{n\rightarrow\infty}\frac{\ln\left(S_{\mu}^{q}(\eta_{n})\right)}{\ln(\eta_{n})}\\
 & \geq &
 \frac{1}{1+R}\\
 & > & \frac{1}{2}.
 \end{eqnarray*}

\end{proof}

\section{R\'{e}nyi dimension of convolutions}

Barbaroux, Germinet and Tcheremchantsev (\cite{BarGerTch}) establish
bounds that relate $D_{q}^{\pm}(\mu\ast\nu)$ with $D_{q}^{\pm}(\mu)$
and $D_{q}^{\pm}(\nu),$ when $q$ is positive and $q\neq1.$ Here
we establish related bounds.

\begin{thm}
Suppose $\mu$ and $\nu$ are Borel measures on $\mathbb{R}^{d}.$
If $1<q<\infty$ and $\frac{1}{q}+1=\frac{1}{r}+\frac{1}{s}$ for
some $1<r,s<\infty,$ then 
\[
D_{q}^{-}(\mu\ast\nu)
\geq\frac{q(r-1)}{r(q-1)}D_{r}^{-}(\mu)+\frac{q(s-1)}{s(q-1)}D_{s}^{-}(\nu).
\]
If $0<q<1,$ and $\mu$ is $q$=finite, and $\frac{1}{q}+1=\frac{1}{r}+\frac{1}{s}$
for some $0<r,s<1,$ then 
\[
D_{q}^{+}(\mu\ast\nu)
\leq\frac{q(r-1)}{r(q-1)}D_{r}^{+}(\mu)+\frac{q(s-1)}{s(q-1)}D_{s}^{+}(\nu).
\]

\end{thm}

\begin{proof}
We let $g$ again be the standard Gaussian. 

Assume $1<r,s<q.$ By Young's inequality (Theorem 1.2.12 in \cite{Grafakos}),
\begin{eqnarray}
 &  & d+\frac{q}{q-1}
 \frac{\ln(\| g_{\sqrt{2}\epsilon}\ast(\mu\ast\nu)\|_{p})}{\ln(\sqrt{2}\epsilon)
 -\ln(\sqrt{2})}\nonumber \\
 & = & d+\frac{q}{q-1}
 \frac{\ln(\|(g_{\epsilon}\ast\mu)\ast(g_{\epsilon}\ast\nu)\|_{p})}{\ln(\epsilon)}
 \nonumber \\
 & \geq & d+\frac{q}{q-1}\left(\frac{\ln(\| g_{\epsilon}\ast\mu\|_{r})}{\ln(\epsilon)}
 +\frac{\ln(\| g_{\epsilon}\ast\nu\|_{s})}{\ln(\epsilon)}\right)
 \label{eq:convComb}\\
 & = & \frac{q(r-1)}{r(q-1)}
 \left(d+\frac{r}{r-1}\frac{\ln(\| g_{\epsilon}\ast\mu\|_{r})}{\ln(\epsilon)}\right).
 \nonumber \\
 &  & \quad+\frac{q(s-1)}{s(q-1)}
 \left(d+\frac{s}{s-1}\frac{\ln(\| g_{\epsilon}\ast\nu\|_{s})}{\ln(\epsilon)}\right)
 \nonumber 
 \end{eqnarray}
Notice that 
\[
0<\frac{q(r-1)}{r(q-1)},\frac{q(s-1)}{s(q-1)}
\]
and
\[
\frac{q(r-1)}{r(q-1)}+\frac{q(s-1)}{s(q-1)}=1.
\]

Now take $\liminf$ of both sides of (\ref{eq:convComb}) and apply
Gu\'{e}rin's formula (\ref{eq:  Guerin's formula}).

For $0<q<1,$ the inequality switches in (\ref{eq:convComb}), we
apply $\limsup$ to both sides, and then use Lemma \ref{lemma: norm of gaussians}.
\end{proof}

\section{Best fit slopes \label{sec:Best-Fit-Slopes}}

Rather than tracking long-term slopes
\[
\frac{\ln\left(S_{\mu}^{q}(e^{x})\right)
-\ln\left(S_{\mu}^{q}(e^{x_{0}})\right)}{x-x_{0}}
\]
to determine a ``dimension'' for $\mu,$ one could look at slope information
of the function
\[
x\mapsto\ln\left(S_{\mu}^{q}(e^{x})\right).
\]
in many ways. Here we consider the slopes of least-squares best fit
lines.

This is not to advocate for or against using least-squares best fit
lines for calculating fractal dimensions in practice.
See \cite{KenkelWalker,KenkelWalkerHTML}
for a discussion. See also \cite{Mandelbrot}.

\begin{rem}
Given a measurable function
\[
\rho:[0,\infty)\rightarrow\mathbb{R}
\]
that is bounded on bounded intervals, the slope of closest line over
$[0,x]$ to $\rho$ is
\[
\frac{6}{x^{3}}\int_{0}^{x}(2t-x)\rho(t)\, dt.
\]
Here closest means with respect to the $L^{2}$ norm using Lebesgue
measure on $[0,x].$

\end{rem}
~

\begin{rem}
Given a sequence
\[
\rho:\mathbb{N}\rightarrow\mathbb{R},
\]
and any positive real $\lambda,$ the slope of the least-squares best
fit line to
\[
\left\{ (0,\rho_{0}),(\lambda,\rho_{1}),\ldots,((n-1)\lambda,\rho_{n-1})\right\} 
\]
is
\[
\frac{6}{(n^{3}-n)\lambda}\sum_{k=0}^{n-1}(2k+1-n)\rho_{k}
\]

\end{rem}
The following lemmas are needed to work for the partition function,
but only depend on bounds from Lemma \ref{thm: order 1 diff}.

\begin{defn}
A function $\rho:[0,\infty)\rightarrow\mathbb{R}$ is \emph{
nearly Lipschitz} if $\rho$ is measurable and there are finite constants
$A$ and $B$ so that
\[
|\rho(x)-\rho(y)|\leq A+B|x-y|.
\]

\end{defn}

\begin{lem}
If $\lambda>0$ and $\rho:[0,\infty)\rightarrow\mathbb{R}$ is 
nearly Lipschitz, then the sequence
\[
n\mapsto\frac{6}{(n\lambda)^{3}}\int_{0}^{n\lambda}(2t-n\lambda)\rho(t)\, dt
\]
 has the same accumulation points as the net
 
\[
x\mapsto\frac{6}{x^{3}}\int_{0}^{x}(2t-x)\rho(t)\, dt.
\]
Moreover, there is a constant $C$ so that
\[
1\leq x\leq y\leq x+\lambda
\]
implies
\[
\left|\frac{6}{x^{3}}\int_{0}^{x}(2t-x)\rho(t)\, dt
-\frac{6}{y^{3}}\int_{0}^{y}(2t-y)\rho(t)\, dt\right|\leq\frac{C}{x}.
\]

\end{lem}

\begin{proof}
Let $A_{0}=A+|\rho(0)|,$ so that
\[
|\rho(t)|\leq A_{0}+Bt.
\]
A change of variable shows
\[
\frac{6}{x^{3}}\int_{0}^{x}(2t-x)\rho(t)\, dt
=6\int_{0}^{1}(2t-1)\frac{1}{x}\rho(xt)\, dt.
\]
For $x>0$ and $y>0,$

\begin{eqnarray*}
 &  & \left|\int_{0}^{1}(2t-1)\frac{1}{x}\rho(xt)\, dt
 -\int_{0}^{1}(2t-1)\frac{1}{y}\rho(yt)\, dt\right|\\
 & = & \left|\int_{0}^{1}(2t-1)\frac{1}{x}\left(\rho(xt)-\rho(yt)\right)\, dt
 +\int_{0}^{1}\left(\frac{1}{x}-\frac{1}{y}\right)(2t-1)\rho(yt)\, dt\right|\\
 & \leq & \frac{1}{x}\int_{0}^{1}|2t-1||\rho(xt)-\rho(yt)|\, dt
 +\frac{|x-y|}{xy}\int_{0}^{1}|2t-1||\rho(yt)|\, dt.
\end{eqnarray*}
We estimate the first term via
\begin{eqnarray*}
\int_{0}^{1}|2t-1||\rho(xt)-\rho(yt)|\, dt 
& \leq & \int_{0}^{1}|\rho(xt)-\rho(yt)|\, dt\\
 & \leq & \int_{0}^{1}A+B|x-y|t\, dt\\
 & = & A+\frac{1}{2}B|x-y|.
\end{eqnarray*}
For the second,
\begin{eqnarray*}
\int_{0}^{1}|2t-1||\rho(yt)|\, dt & \leq & \int_{0}^{1}|\rho(yt)|\, dt\\
 & \leq & \int_{0}^{1}A_{0}+Byt\, dt
\end{eqnarray*}
so
\begin{equation}
\int_{0}^{1}|2t-1||\rho(yt)|\, dt\leq A_{0}+\frac{1}{2}By.
\label{eq:overall integral bound}
\end{equation}
Therefore
\begin{eqnarray*}
 &  & \left|\int_{0}^{1}(2t-1)\frac{1}{x}\rho(xt)\, dt
 -\int_{0}^{1}(2t-1)\frac{1}{y}\rho(yt)\, dt\right|\\
 & \leq & \frac{1}{x}\left(A+\frac{1}{2}B|x-y|\right)
 +\frac{|x-y|}{xy}\left(A_{0}+\frac{1}{2}By\right).
\end{eqnarray*}
If we assume
\[
x\leq y\leq x+\lambda
\]
then
\begin{eqnarray*}
 &  & \left|\int_{0}^{1}(2t-1)\frac{1}{x}\rho(xt)\, dt
 -\int_{0}^{1}(2t-1)\frac{1}{y}\rho(yt)\, dt\right|\\
 & \leq & \frac{1}{x}\left(A+\frac{1}{2}B\lambda\right)
 +\frac{\lambda}{xy}\left(A_{0}+\frac{1}{2}By\right)\\
 & \leq & \frac{1}{x}\left(A+\frac{1}{2}B\lambda\right)
 +\frac{\lambda}{x^{2}}A_{0}+\frac{1}{2}\frac{\lambda}{x}B\\
 & = & \left(A+B\lambda\right)\frac{1}{x}
 +\left(\lambda A_{0}\right)\frac{1}{x^{2}}
\end{eqnarray*}

\end{proof}

\begin{lem}
If $\lambda>0$ and $\rho:[0,\infty)\rightarrow\mathbb{R}$ is 
nearly Lipschitz, then the sequence
\[
n\mapsto\frac{6}{(n^{3}-n)\lambda}\sum_{k=0}^{n-1}(2k+1-n)\rho(\lambda k)
\]
 has the same accumulation points as the sequence
 
\[
n\mapsto\frac{6}{(n\lambda)^{3}}\int_{0}^{n\lambda}(2t-n\lambda)\rho(t)\, dt.
\]
Moreover, there is a constant $C$ so that
\[
2\leq n
\]
implies
\[
\left|\frac{6}{(n\lambda)^{3}}\int_{0}^{n\lambda}(2t-n\lambda)\rho(t)\, dt
-\frac{6}{(n^{3}-n)\lambda}\sum_{k=0}^{n-1}(2k+1-n)\rho(\lambda k)\right|
\leq\frac{C}{n}.
\]

\end{lem}

\begin{proof}
Let $A_{0}=A+|\rho(0)|,$ so that
\[
|\rho(t)|\leq A_{0}+Bt.
\]

Equation (\ref{eq:overall integral bound}), with $y=n\lambda,$ gives
us the constant bound
\begin{eqnarray*}
\left|\frac{6}{(n\lambda)^{3}}\int_{0}^{n\lambda}(2t-n\lambda)\rho(t)\, dt\right| 
& = & \frac{6}{n\lambda}\left|\int_{0}^{1}(2t-1)\rho(n\lambda t)\, dt\right|\\
 & \leq & \frac{6}{n\lambda}\left(A_{0}+\frac{1}{2}Bn\lambda\right)\\
 & \leq & 3\left(\frac{A_{0}}{\lambda}+B\right)
 \end{eqnarray*}
as long as $n\geq2.$ Also
\[
2\leq n\implies\left|\frac{n^{3}}{n^{3}-n}-1\right|\leq\frac{2}{n^{2}}.
\]
Therefore, it suffices to estimate the distance from
\[
\frac{6}{(n\lambda)^{3}}\int_{0}^{n\lambda}(2t-n\lambda)\rho(t)\, dt
=\frac{6}{n^{3}\lambda}\int_{0}^{n}(2t-n)\rho(\lambda t)\, dt
\]
to
\[
\frac{6}{n^{3}\lambda}\sum_{k=0}^{n-1}(2k+1-n)\rho(\lambda k).
\]

For all $n\geq2,$
\begin{eqnarray*}
 &  & \left|\frac{6}{n^{3}\lambda}\int_{0}^{n}(2t-n)\rho(\lambda t)\, dt
 -\frac{6}{n^{3}\lambda}\sum_{k=1}^{n-1}(2k+1-n)\rho(\lambda k)\right|\\
 & = & \frac{6}{n^{3}\lambda}
 \left|\sum_{k=0}^{n-1}\int_{k}^{k+1}(2t-n)\rho(\lambda t)\, dt
 -\sum_{k=0}^{n-1}\int_{k}^{k+1}(2k+1-n)\rho(\lambda k)\, dt\right|\\
 & \leq & \frac{6}{n^{3}\lambda}
 \sum_{k=0}^{n-1}\int_{k}^{k+1}|(2t-n)-(2k+1-n)||\rho(\lambda t)|\, dt\\
 &  & \quad+\frac{6}{n^{3}
 \lambda}\sum_{k=0}^{n-1}
 \int_{k}^{k+1}|2k+1-n||\rho(\lambda t)-\rho(\lambda k)|\, dt\\
 & \leq & \frac{6}{n^{3}\lambda}
 \sum_{k=0}^{n-1}
 \int_{k}^{k+1}|\rho(\lambda t)|+n|\rho(\lambda t)-\rho(\lambda k)|\, dt\\
 & \leq & \frac{6}{n^{3}\lambda}
 \sum_{k=0}^{n-1}\int_{k}^{k+1}A_{0}+B\lambda t+n(A+B\lambda(t-k))\, dt\\
 & = & \left(\frac{6A}{\lambda}
 +\frac{6|\rho(0)|}{\lambda}\right)\frac{1}{n^{2}}+
 \left(9B+\frac{6A}{\lambda}\right)\frac{1}{n}.
 \end{eqnarray*}

\end{proof}

\begin{lem}
\label{lemma:  discrete fit ok}If $\lambda>0$ and 
$\rho:[0,\infty)\rightarrow\mathbb{R}$
is  nearly Lipschitz, then the sequence
\[
n\mapsto\frac{6}{(n^{3}-n)\lambda}\sum_{k=0}^{n-1}(2k+1-n)\rho(\lambda k)
\]
 has the same accumulation points as the net
 
\[
x\mapsto\frac{6}{x^{3}}\int_{0}^{x}(2t-x)\rho(t)\, dt.
\]
Moreover, there are constants $C$ and $D$ so that
\[
D\leq n\lambda\leq x\leq(n+1)\lambda
\]
implies
\[
\left|\frac{6}{x^{3}}\int_{0}^{x}(2t-x)\rho(t)\, dt
-\frac{6}{(n^{3}-n)\lambda}\sum_{k=0}^{n-1}(2k+1-n)\rho(\lambda k)\right|
\leq\frac{C}{n}.
\]

\end{lem}

\begin{proof}
Let $C_{1}$ and $C_{2}$ be the constants from the last two lemmas.
Suppose
\[
\max(1,2\lambda)\leq n\lambda\leq x\leq(n+1)\lambda.
\]
Let $y=n\lambda.$ We have
\[
2\leq n
\]
and
\[
0<y\leq x\leq y+\lambda.
\]
Therefore
\begin{eqnarray*}
 &  & \left|\frac{6}{x^{3}}\int_{0}^{x}(2t-x)\rho(t)\, dt
 -\frac{6}{(n^{3}-n)\lambda}\sum_{k=0}^{n-1}(2k+1-n)\rho(\lambda k)\right|\\
 & \leq & \left|\frac{6}{x^{3}}\int_{0}^{x}(2t-x)\rho(t)\, dt
 -\frac{6}{y^{3}}\int_{0}^{y}(2t-y)\rho(t)\, dt\right|\\
 &  &\quad +\left|\frac{6}{n^{3}\lambda^{3}}\int_{0}^{n\lambda}(2t-n\lambda)\rho(t)\, dt
 -\frac{6}{(n^{3}-n)\lambda}\sum_{k=0}^{n-1}(2k+1-n)\rho(\lambda k)\right|\\
 & \leq & \frac{C_{1}}{y}+\frac{C_{2}}{n}\\
 & = & \left(\frac{C_{1}}{\lambda}+C_{2}\right)\frac{1}{n}.
\end{eqnarray*}

\end{proof}

\begin{thm}
Suppose $\mu$ is a finite Borel measure on $\mathbb{R}^{d},$ and
$0<q<1$ or $1<q<\infty.$ If $0<q<1$ then also suppose $\mu$ is
$q$-finite. Suppose $v>1.$ Let 
\[
S_{\mu}^{q}(\epsilon)
=\sum_{k\in\mathbb{Z}}\mu(\epsilon k+\epsilon\mathbb{I})^{q}.
\]
For $x$ in $[0,\infty),$ let $m_{x}$ and $b_{x}$ be the real numbers
so that
\[
t\mapsto m_{x}t+b_{x}\]
is the linear function that minimizes
\[
\int_{-x}^{0}
\left(\ln\left(S_{\mu}^{p}(e^{t})\right)-\left(m_{x}t+b_{x}\right)\right)^{2}\: dt.
\]
(Here $dt$ refers to Lebesgue measure on $[-x,0].$) For a natural
number $n,$ let $\tilde{m}_{n}$ and $\tilde{b}_{n}$ the the real
numbers so that
\[
t\mapsto\tilde{m}_{n}t+\tilde{b}_{n}\]
is the the least-squares best fit line to the pairs 
\[
\left\{ \left.\left(
\ln(v^{-k}),\ln\left(S_{\mu}^{q}(v^{-k})\right)\right)\right|k=0,1,\ldots,n-1
\right\} .
\]
There are constants $C$ and $D$ so that
\[
D\leq n\ln(v)\leq x\leq(n+1)\ln(v)
\]
implies
\[
\left|m_{x}-\tilde{m}_{n}\right|\leq\frac{C}{n}.
\]
In particular,
\[
\lim_{x\rightarrow\infty}\,_{\textrm{inf }}^{\textrm{sup}}m_{x}
=\lim_{n\rightarrow\infty}\,_{\textrm{inf }}^{\textrm{sup}}\tilde{m}_{n}.
\]

\end{thm}

\begin{proof}
By Theorem \ref{thm: order 1 diff}, there are constants $A$ and
$B$ so that
\[
\left|\ln\left(S_{\mu}^{q}(e^{-x})\right)
-\ln\left(S_{\mu}^{q}(e^{-y})\right)\right|\leq A+B\left|x-y\right|
\]
for all $x$ and $y.$ Since $\ln\left(S_{\mu}^{q}(e^{-x})\right)$
is Borel measurable, we see from Lemma \ref{lemma:  discrete fit ok}
that there is are constants $C$ and $D$ so that
\[
D\leq n\ln(\nu)\leq x\leq(n+1)\ln(\nu)
\]
implies
\[
\frac{6}{x^{3}}\int_{0}^{x}
(2t-n\lambda)\ln\left(S_{\mu}^{p}\left(e^{-t}\right)\right)\, dt
\]
is within $\frac{C}{n}$ of
\[
\frac{6}{(n^{3}-n)\ln(\nu)}\sum_{k=0}^{n-1}
(2k+1-n)\ln\left(S_{\mu}^{p}\left(e^{-k\ln(\nu)}\right)\right).
\]
That is,\begin{equation}
\frac{6}{x^{3}}\int_{0}^{x}
(2t-n\lambda)\ln\left(S_{\mu}^{p}\left(e^{-t}\right)\right)\, dt
\label{eq:S integral}
\end{equation}
is within $\frac{C}{n}$ of
\begin{equation}
\frac{6}{(n^{3}-n)\ln(\nu)}\sum_{k=0}^{n-1}
(2k+1-n)\ln\left(S_{\mu}^{p}\left(\nu^{-k}\right)\right).
\label{eq: S sum}
\end{equation}
The quantity in (\ref{eq:S integral}) gives the slope of the best
fit over $[0,x]$ of
\[
t\mapsto\ln\left(S_{\mu}^{q}(e^{-t})\right),
\]
and so\begin{equation}
m_{x}=\frac{-6}{x^{3}\ln(\nu)}\int_{0}^{x}(2t-x)\ln\left(S_{\mu}^{q}(e^{-t})\right)\, dt.
\label{eq:mx one}
\end{equation}
The quantity in (\ref{eq: S sum}) gives the slope of the best fit
to
\[
\left\{ \left.\left(
k\ln(v),\ln\left(S_{\mu}^{q}(v^{-k})\right)\right)\right|k=0,1,\ldots,n-1
\right\} 
\]
 and so
 \begin{equation}
\tilde{m}_{n}=\frac{-6}{(n^{3}-n)\ln(\nu)}
\sum_{k=0}^{n-1}(2k+1-n)\ln\left(S_{\mu}^{q}(\nu^{-k})\right).\label{eq: m tilde x 1}
\end{equation}
We are done.
\end{proof}

\begin{rem}
It is interesting to note some alternative formulas:
\begin{equation}
\tilde{m}_{n}=\frac{\sum_{k=0}^{n-1}(2k+1-n)
\ln\left(S_{\mu}^{q}(\nu^{-k})\right)}{\sum_{k=0}^{n-1}(2k+1-n)\ln\left(\nu^{-k}\right)};
\label{eq: m tilde x 2}
\end{equation}
\begin{equation}
\tilde{m}_{n}=\frac
{6}{(n^{3}-n)\ln(\nu)}
\sum_{k=1}^{n-1}k(n-k)\left(\ln\left(S_{\mu}^{q}(\nu^{-(k-1)})
-\ln\left(S_{\mu}^{q}(\nu^{-k})\right)\right)\right);
\label{eq:m tilde x 3}
\end{equation}
\begin{equation}
\tilde{m}_{n}=\frac
{\sum_{k=1}^{n-1}k(n-k)\left(\ln\left(S_{\mu}^{q}(\nu^{-k})\right)
-\ln\left(S_{\mu}^{q}(\nu^{-(k-1)})\right)\right)}
{\sum_{k=1}^{n-1}k(n-k)\ln\left(\nu^{-1}\right)};\label{eq: m tilde x 4}\end{equation}
\begin{equation}
m_{x}=\frac
{\int_{0}^{x}(2t-(x))\ln\left(S_{\mu}^{q}(\nu^{-t})\right)\, dt}
{\int_{0}^{x}(2t-x)\ln\left(\nu^{-t}\right)\, dt}.\label{eq:mx 2}
\end{equation}

\end{rem}

\section{More examples \label{sec:Best-Fit-Slopes II}}

The slope of the least-squares best fit linear approximation
to the partition function cannot always be used to determine
the Renyi dimension.  We show this by the following example.
This example is far from what we hope to see in applications.

\begin{lem}
For any $1<q<\infty,$ there is a finite Borel measure $\mu$ on $[0,1]$
for which 
\[
\limsup_{x\rightarrow-\infty}
\frac{\ln\left(S_{\mu}^{q}(e^{x})\right)}{x}<\limsup_{x\rightarrow-\infty}m_{x},
\]
where
\[
m_{x}t+b_{x}\approx\ln\left(S_{\mu}^{q}(e^{t})\right)\quad(-x\leq t\leq0)
\]
 is the least-squares best fit line. More specifically, $m_{x}$ and
$b_{x}$ minimize
\[
\int_{-x}^{0}\left(\ln\left(S_{\mu}^{q}(e^{t})\right)-(m_{x}t+b_{x})\right)^{2}\, dt.
\]
For any $0<q<1,$ there is a finite Borel measure $\mu$ on $[0,1]$
for which 
\[
\liminf_{\epsilon\rightarrow0}
\frac{\ln\left(S_{\mu}^{q}(e^{x})\right)}{x}>\liminf_{x\rightarrow-\infty}m_{x},
\]
\end{lem}

\begin{proof}
We will use the $\mu$ associated with a sequence $a_{n},$ as in
Lemma \ref{lem:a_n example}.

Let 
\[
a_{1}=\frac{30}{47}
\]
and
\[
a_{k}=\left\{ \begin{array}{cc}
0 
 & \textrm{if }48^{n}<k\leq12(48^{n}),\textrm{ any }n\in\mathbb{N}\\
1 
 & \textrm{if }12(48^{n})<k\leq36(48^{n}),\textrm{ any }n\in\mathbb{N}\\
\frac{1}{2} 
 & \textrm{if }36(48^{n})<k\leq48(48^{n}),\textrm{ any }n\in\mathbb{N}\end{array}\right..
\]
When $1<q<\infty,$
\[
\limsup_{\epsilon\rightarrow0}\frac{\ln\left(S_{\mu}^{q}(\epsilon)\right)}{\ln(\epsilon)}
=\limsup_{n\rightarrow0}(q-1)\ln(2)\frac{1}{n}\sum_{k=1}^{n}a_{k}
\]
and
\[
\limsup_{\epsilon\rightarrow0}m_{x}=\limsup_{n\rightarrow0}\frac{6(q-1)\ln(2)}{n^{3}}\sum_{k=1}^{n-1}k(n-k)a_{k}.
\]
Therefore we really only need to show that
\[
\limsup_{n\rightarrow\infty}\frac{1}{n}\sum_{k=1}^{n}a_{k}<\limsup_{n\rightarrow\infty}\frac{6}{n^{3}}\sum_{k=1}^{n-1}k(n-k)a_{k}.
\]
For $q<1$ the $(q-1)$ reverses the inequalities and turns each $\liminf$
into a $\limsup,$ so the desired inequality reduces to the same thing.

For all $n,$
\begin{eqnarray*}
\frac{1}{48^{n+1}}\sum_{k=1}^{48^{n+1}}a_{k} 
& = & \frac{1}{48^{n+1}}
\left[\sum_{k=1}^{48^{n}}a_{k}+24(48^{n})+12(48^{n})\frac{1}{2}\right]\\
 & = & \frac{1}{48}\left[\frac{1}{48^{n}}\sum_{k=1}^{48^{n}}a_{k}+30\right].
 \end{eqnarray*}
Since
\[
\frac{1}{48}\left[\frac{30}{47}+30\right]=\frac{30}{47}
\]
we have found
\[
\frac{1}{48^{n}}\sum_{k=1}^{48^{n}}a_{k}=\frac{30}{47}
\]
 for all $n.$ The terms of value $0$ will cause the average to fall
until index $12(48^{n}).$ At this point, the average will be
\[
\frac{1}{12}\frac{30}{47}+\frac{11}{12}0=\frac{5}{94}
\]
The next $12(48^{n})$ terms are of value $1,$ so the average rises
to
\[
\frac{1}{3}\frac{5}{94}+\frac{2}{3}1=\frac{193}{282}.
\]
Next the average falls, due the terms of value $\frac{1}{2},$ until
it is back to $\frac{30}{47}.$ Therefore,
\[
\limsup_{n\rightarrow\infty}\frac{1}{n}\sum_{k=1}^{n}a_{k}=\frac{193}{282}.
\]

We now need just a decent estimate on
\[
\limsup_{n\rightarrow\infty}\frac{6(p-1)}{n^{3}}\sum_{k=1}^{n-1}k(n-k)a_{k}.
\]
Indeed, since $a_{n}>0$ for all $n,$ we find
\begin{eqnarray*}
 &  & \limsup_{n\rightarrow\infty}\frac{6}{n^{3}}\sum_{k=1}^{n-1}k(n-k)a_{k}\\
 & \geq & \limsup_{n\rightarrow\infty}\frac{6}{(48^{n})^{3}}
 \sum_{k=1}^{48^{n}}k(48^{n}-k)a_{k}\\
 & > & \limsup_{n\rightarrow\infty}
 6\frac{\sum_{k=12(48^{n-1})+1}^{36(48^{n-1})}k(48^{n}-k)
 +\sum_{k=36(48^{n-1})+1}^{48^{n}}k(48^{n}-k)\frac{1}{2}}{(48^{n})^{3}}\\
 & = & \lim_{n\rightarrow\infty}
 6\sum_{k=12(48^{n-1})+1}^{36(48^{n-1})}
 \frac{k}{48^{n}}\left(1-\frac{k}{48^{n}}\right)\frac{1}{48^{n}}\\
 &  & \quad+\lim_{n\rightarrow\infty}3\sum_{k=36(48^{n-1})+1}^{48^{n}}
 \frac{k}{48^{n}}\left(1-\frac{k}{48^{n}}\right)\frac{1}{48^{n}}\\
 & = & 6\int_{\frac{1}{4}}^{\frac{3}{4}}t\left(1-t\right)\, dt
 +3\int_{\frac{3}{4}}^{1}t\left(1-t\right)\, dt\\
 & = & 6\left(\frac{1}{2}\left(\frac{3}{4}\right)^{2}
 -\frac{1}{3}\left(\frac{3}{4}\right)^{3}\right)
 -3\left(\frac{1}{2}\left(\frac{1}{4}\right)^{2}
 -\frac{1}{3}\left(\frac{1}{4}\right)^{3}\right)\\
 & = & \frac{49}{64}.
 \end{eqnarray*}
Thus
\[
\limsup_{\epsilon\rightarrow0}\frac{\ln\left(S_{\mu}^{q}(\epsilon)\right)}{\ln(\epsilon)}
=\frac{193}{282}(q-1)=\frac{6176}{9024}(q-1)
\]
and
\[
\limsup_{\epsilon\rightarrow0}m_{x}>\frac{49}{64}(q-1)
=\frac{6909}{9024}(q-1).
\]
 
\end{proof}

\section{Modified R\'{e}nyi Dimensions}

The theory of regular variation and its extensions (\cite{BingGoldTeu})
give many ways to measure how closely a function $f$ behaves like
various powers $x^{c}$ near $\infty.$ Regular variation forces $f$
to behave like a single power $x^{\delta}.$ 

More realistic classes are those of extended variation and $O$-regularly
varying functions. Both classes allow $f$ to behave like $x^{c}$
for $c$ in a range $(\alpha,\beta),$ but they differ on the meaning
of ``behave.'' (This is a bit vague. See\cite{BingGoldTeu}.) The
extended real numbers $\alpha$ and $\beta$ are called Karamata indices
in the case where $f$ is of extended variation. For the class of
$O$-regularly varying functions, these are called the Matuszewska
indices.

Guido and Isola (\cite{GuidoIsola1,GuidoIsola2,GuidoIsola3}) have
used the Matuszewska indices to define a new local fractal dimension.
Stern (\cite{Stern}) has suggested generally that the theories of
extended variation and $O$-regular variation be applied to global
fractal dimensions.

The example in Section \ref{sec:Best-Fit-Slopes II} is rather unnatural.
It can perhaps be explained away if we use Matuszewska indices to
describe the ``slope at infinity'' of the partition function.

We use the following as a working definition of the Matuszewska indices.
It is equivalent to the standard definition, c.f. pages 68--73 of
\cite{BingGoldTeu}.

\begin{defn}
Suppose
\[
f:[0,\infty)\rightarrow(0,\infty)
\]
is any function. The \emph{upper Matuszewska index of} $f$ \emph{is
\[
\alpha(f)=\left\{ \alpha\in\mathbb{R}\left|\exists X,C\textrm{ s.t. }y
\geq x
geq X\implies f(y)\leq f(x)C\left(\frac{y}{x}\right)^{\alpha}\right.\right\} .
\]
Here $X$ are $C$ are to be understood to be positive real numbers.}
The \emph{lower Matuszewska index of} $f$ \emph{is
\[
\beta(f)=\left\{ \beta\in\mathbb{R}\left|\exists X,C\textrm{ s.t. }y
\geq x
\geq X\implies f(y)\geq f(x)C\left(\frac{y}{x}\right)^{\beta}\right.\right\} .
\]
}
\end{defn}

It is easy to show that
\[
\beta(f)\leq\liminf_{x\rightarrow\infty}\frac{\ln(f(x))}{\ln(x)}
\leq\limsup_{x\rightarrow\infty}\frac{\ln(f(x))}{\ln(x)}\leq\alpha(f).
\]
Again, see \cite{BingGoldTeu}. The middle numbers are the so-called
order of $f.$ In an unfortunate clash of terminology, the ``upper
and lower R\'{e}nyi dimensions of order $q$'' are the upper and
lower orders of 
\[
x\mapsto\left(S_{\mu}^{q}(x^{-1})\right)^{\frac{1}{1-q}}.
\]
Perhaps it is better to refer to $q$ as the index.

\begin{defn}
If $\mu$ is a finite measure, and if $0<q<1$ or $1<q<\infty,$ the
\emph{upper} and \emph{lower Matuszewska Dimensions of index} $q$
are the upper and lower Matuszewska indices of the function
\[
x\mapsto\left(S_{\mu}^{q}(x^{-1})\right)^{\frac{1}{1-q}},
\]
denoted $D_{q}^{++}(\mu)$ and $D_{q}^{--}(\mu)$ respectively.
\end{defn}

\begin{thm}
Suppose $\mu$ is a finite Borel measure on $\mathbb{R}^{d}.$ If
$1<q<\infty,$ or if $0<q<1$ and $\mu$ is $q$-finite, then the
partition function $S_{\mu}^{q}(x^{-1})$ is of extended variation
and
\[
0\leq D_{q}^{--}(\mu)\leq D_{q}^{-}(\mu)\leq D_{q}^{+}(\mu)\leq D_{q}^{++}(\mu)\leq d.
\]

\end{thm}

\begin{proof}
The second and fourth inequalities come from the general facts about
order and Matuszewska indices. The middle is even more standard. The
outer inequalities are really just restatements of those in Theorem
\ref{thm: order 1 diff}.

Equivalently, these inequalities show that the upper and lower Matuszewska
indices of $S_{\mu}^{q}(x^{-1})$ are bounded between and $0$ and
$1-q.$ Since $S_{\mu}^{q}(x^{-1})$ is measurable, we can apply
\cite[Theorem 2.1.7]{BingGoldTeu}
to conclude that $S_{\mu}^{q}(x^{-1})$ is of extended variation.
\end{proof}

\begin{rem}
In the example of Section \ref{sec:Best-Fit-Slopes II}:
\[
D_{q}^{--}(\mu)=0,
\]
\[
D_{q}^{-}(\mu)=\frac{5}{94},
\]
\[
D_{q}^{+}(\mu)=\frac{193}{282},
\]
\[
D_{q}^{++}(\mu)=1.
\]
Thus the upper and lower Matuszewska dimensions dismiss $\mu$ from
``being fractal'' more resoundingly than do the upper and lower R\'{e}nyi
dimensions.
\end{rem}

\end{document}